\documentstyle[12pt]{article}
\catcode`\@=11
\@addtoreset{equation}{section}
\def\theequation{\thesection.\arabic{equation}}
\catcode`\@=12
\newtheorem{Theorem}{Theorem}[section]
\newtheorem{Definition}{Definition}[section]
\newtheorem{Proposition}{Proposition}[section]
\newtheorem{Lemma}{Lemma}[section]
\newtheorem{Corollary}{Corollary}[section]

\title{A Proof On  
Weinstein Conjecture On Cotangent Bundles
\thanks{Project 19871044 Supported by NSF}}

\author{Renyi Ma \\
Department of Mathematics \\
Tsinghua University \\
Beijing, 100084\\
People's Republic of China\\
rma@math.tsinghua.edu.cn}

\date { }

\begin{document}
\textwidth=165mm
\textheight=200mm
\parindent=8mm
\frenchspacing
\maketitle

\begin{abstract}
In this article, we prove that there 
exists at least one closed characteristics 
of Reeb vector field in a connected contact manifolds of induced type in 
the cotangent bundles of any open smooth manifolds which 
confirms completely the Weinstein conjecture in cotangent 
bundles of open manifold. 
\end{abstract}


\noindent{\bf 2000 MR Subject Classification} 32Q65, 53D35, 53D12

\section{Introduction and results }

A contact structure on a manifold  is a field of a tangent
hyperplanes (contact hyperplanes) that is nondegenerate at any point.
Locally such a field is defined as the field of zeros of a $1-$form
$\lambda $, called a contact form. The nondegeneracy condition
is that $d\lambda $ is nondegenerate on the hyperplanes on which $\lambda $
vanishes; equivalently, in $(2n-1)-$space:

$$\lambda \wedge (d\lambda )^{n-1}\neq 0$$

The important example of contact manifold is
the well-known projective cotangent bundles definded as
follows:

Let $N=T^*M$ be the cotangent bundle of a smooth connected compact
manifold $M$. $N$ carries a canonical symplectic structure $\omega
=-d\lambda $ where $\lambda=\sum _{i=1}^{n}y_idx_i$ is the Liouville
form on $N$, see \cite{ag,kl}. Let $P=PT^*M$ be the oriented projective
cotangent bundle of $M$, i.e. $P=\cup _{x\in M}PT^*_xM$. It is well
known that $P$ carries a canonical contact structure induced by the
Liouville form and the projection $\pi :T^*M\mapsto PT^*M$.

Let $(\Sigma ,\lambda )$ 
be a smooth closed oriented manifold of dimension
$2n-1$ with a contact form $\lambda $.
Associated to $\lambda$ there are two important structures.
First of all the so-called Reed vectorfield $X_\lambda $ defined
by
$$i_{X_\lambda }\lambda   \equiv 1, \ \ i_{X_{\lambda }}d\lambda  \equiv 0$$
and secondly the contact structure $\xi =\xi _{\lambda }
\mapsto \Sigma $ given by
$$\xi _{\lambda }=\ker (\lambda )\subset T \Sigma $$
by a result of Gray, \cite{gra} , the contact structure is very
stable. In fact, if $(\lambda  _t  )_{t\in  [0,1]}$ is a smooth
arc of contact forms inducing the arc of contact structures
$(\xi _t)_{t\in [0,1]}$, there exists a smooth arc
$(\psi _t)_{t\in [0,1]}$
of diffeomorphisms with $\psi _0=Id$, such that
\begin{equation}
T\Psi _t(\xi _0)=\xi _t
\end{equation}
here it is important that $\Sigma $ is compact. From (1.1) and
the fact that $\Psi _0=Id$ it follows immediately that there
exists a smooth family of maps $[0,1]\times \Sigma \mapsto
(0,\infty ):(t, m)\to f_t(m)$ such that
\begin{equation}
\Psi ^*_t\lambda _t=f_t\lambda _0
\end{equation}
In contrast to the contact structure the dynamics of the
Reeb vectorfield changes drastically under small perturbation
and in general the flows associated to $X_t$ and $X_s$ for
$t\neq s$ will not be conjugated, see\cite{ag,eg}.

\vskip 3pt 

  Let $M$ be a Riemann manifold with Riemann metric, then it is well known 
that there exists a canonical contact structure in the unit 
sphere of its tangent bundle and the motion of geodesic line lifts 
to a geodesic flow on the unit sphere bundles. Therefore 
the closed orbit of geodesic flow or Reeb flow 
on the sphere bundle projects to a closed geodesics in 
the Riemann manifolds, conversely the closed geodesic orbit lifts 
to a closed Reeb orbit. The classical work of Ljusternik and Fet 
states that every simply connected Riemannian
manifold has at least one closed geodesics, this with the 
Cartan and Hadamard's results on non-simply closed Riemann 
manifold implies that any closed 
Riemann manifolds has a closed geodesics, i.e., 
the sphere bundle of a closed Riemann manifold with 
standard contact form carries at least one closed Reeb orbits 
which is a lift of closed geodesics of base manifold. 
Its proof depends on the classical minimax principle 
of Ljusternik and Schnirelman or minimalization 
of Hadamard and Cartan,\cite{kl}, an 
$J-$holomorphic curve's proof can be found in \cite{ma4}. 
In sympletic geometry, Gromov \cite{gro} 
introduces the global methods to proves the existences of 
symplectic fixed points or periodic orbits which depends on the 
nonlinear Fredholm alternative of $J-$holomorphic curves 
in the symplectic manifolds. In this paper we use the 
$J-$holomorphic curve's method to prove 

\begin{Theorem}
Let $(\Sigma ,\lambda )$ be a contact manifold with 
contact form $\lambda $ of induced type or Weinstein type in the 
cotangent bundles of any open smooth manifold with symplectic 
form $\sum _{i=1}^ndp_i\wedge dq_i$ induced by 
Liouville form $\alpha =\sum _{i=1}^np_idq_i$, i.e., 
there exists a transversal vector field $Z$ to $\Sigma $ 
such that $L_Z\omega =\omega $, $\lambda =i_Z\omega $.   
Let $X_{\lambda } $ be its Reeb vector 
field. Then, 
there exists at least one closed characteristic for 
$X_\lambda $. 
\end{Theorem}

This gives a complete solution on the well-known 
Weinstein conjecture in cotangent bundles of smooth open manifold.
Note that Viterbo \cite{vi1} first proved the 
above result for any contact manifolds $\Sigma $ of 
induced type in $R^{2n}=T^*R^n$ after 
Rabinowitz \cite{ra} and Weinstein \cite{we1,we2}. 
After Viterbo's work many results 
were obtained in \cite{fhv,ho,hv1,hz,ma1,ma2} etc 
by using variational method or Gromov's 
$J-$holomorphic curves via 
nonlinear Fredholm alternative, see survey paper 
\cite{el}. 

\begin{Corollary}
(\cite{ma2})If $M=N\times R$, $N$ is any 
smooth manifold, then Theorem 1.1 holds true.
\end{Corollary}

Through the variational method 
by Hofer and Viterbo\cite{hv1}, especially, Viterbo finally in 
\cite{vi2} proved the following result.
\begin{Corollary}
(Viterbo\cite{vi2})If $M$ is an open $simply$ 
connected manifolds and $[\lambda -p_idq_i]=0$, 
then Theorem 1.1 holds true.
\end{Corollary}

\vskip 3pt 

{\bf Sketch of proofs}: We work in the framework 
as in \cite{gro,ma3}. In Section 2, we study the 
linear Cauchy-Riemann operator and sketch some basic 
properties. In section 3, first we construct a 
Lagrangian submanifold $W$ under the assumption that 
there does not exists closed Reeb orbit in 
$(\Sigma ,\lambda )$; second, we study 
the space ${\cal D}(V,W)$
of contractible disks in manifold $V$ with boundary 
in Lagrangian submanifold $W$ and construct a Fredholm
section of tangent bundle of  ${\cal  D }(V,W)$.  
In section 4, following \cite{gro,ma3}, we prove that the Fredholm 
section is not proper by using a 
special anti-holomorphic section as in \cite{gro,ma3}.  
In section 5, 
we transform the non-homogeneous Cauchy-Riemann equation as 
$J-$holomorphic 
curves.
In the final section, we
use nonlinear Fredholm trick in \cite{gro,ma3}  
to 
complete our proof.

\section{Linear Fredholm theory}

For $100<k<\infty $ consider the Hilbert space
$V_k$ consisting of all maps $u\in H^{k,2}(D, C^n)$,
such that $u(z)\in R^n\subset C^n$ for almost all $z\in
\partial D$. $L_{k-1}$ denotes the usual Hilbert $L_{k-1}-$space
$H_{k-1}(D, C^n)$. We define an operator
$\bar \partial :V_p\mapsto L_p$ by
\begin{equation}
\bar \partial u=u_s+iu_t
\end{equation}
where the coordinates on $D$ are
$(s,t)=s+it$, $D=\{ z||z|\leq 1\}  $.
The following result is well known(see\cite{al,wen}).
\begin{Proposition}
$\bar \partial :V_p\mapsto L_p $ is a surjective real
linear Fredholm operator of index $n$. The kernel
consists of the constant real valued maps.
\end{Proposition}
Let $(C^n, \sigma =-Im(\cdot ,\cdot ))$ be the standard
symplectic space. We consider a real $n-$dimensional plane
$R^n\subset C^n$. It is called Lagrangian if the
skew-scalar product of any two vectors of $R^n$ equals zero.
For example, the plane $p=0$ and $q=0$ are Lagrangian
subspaces. The manifold of all (nonoriented) Lagrangian subspaces of
$R^{2n}$ is called the Lagrangian-Grassmanian $\Lambda (n)$.
One can prove that the fundamental group of
$\Lambda (n)$ is free cyclic, i.e.
$\pi _1(\Lambda (n))=Z$. Next assume
$(\Gamma (z))_{z\in \partial D}$ is a smooth map
associating to a point $z\in \partial D$ a Lagrangian
subspace $\Gamma (z)$ of $C^n$, i.e.
$(\Gamma (z))_{z\in \partial D}$ defines a smooth curve
$\alpha $ in the Lagrangian-Grassmanian manifold $\Lambda (n)$.
Since $\pi _1(\Lambda (n))=Z$, one have
$[\alpha ]=ke$, we call integer $k$ the Maslov index
of curve $\alpha $ and denote it by $m(\Gamma )$, see(\cite{ag}).

Now let $z:S^1\mapsto R^n\subset C^n$ be a smooth curve. 
Then
it defines a constant loop $\alpha $ in Lagrangian-Grassmanian
manifold $\Lambda (n)$. This loop defines
the Maslov index $m(\alpha )$ of the map
$z$ which is easily seen to be zero.

  Now Let $(V,\omega )$ be a symplectic manifold and 
$W\subset V$ a closed Lagrangian submanifold. Let 
$u:D^2\to V$ be a smooth map homotopic to constant map 
with boundary 
$\partial D\subset W$.
Then $u^*TV$ is a symplectic vector bundle and 
$(u|_{\partial D})^*TW$ be a Lagrangian subbundle in 
$(u|_{\partial D})^*TV$. Since $u$ is contractible, we can take 
a trivialization of $u^*TV$ as 
$$\Phi (u^*TV)=D\times C^n$$ 
and 
$$\Phi (u|_{\partial D})^*TW)\subset S^1\times C^n$$
Let 
$$\pi _2: D\times C^n\to C^n$$
then 
$$\bar u: z\in S^1\to 
\{\pi _2\Phi (u|_{\partial D})^*TW(z)\} \in \Lambda (n).$$
Write $\bar u=u|_{\partial D}$.
\begin{Lemma}
Let $u: (D^2,\partial D^2) \rightarrow (V,W)$ be a $C^k-$map $(k\geq 1)$ 
as above. Then,
$$m(u|_{\partial D})=0$$
\end{Lemma}
Proof.  Since $u$ is contractible in $V$ relative to $W$, we have 
a homotopy $\Phi _s$ of trivializations
such that 
$$\Phi _s(u^*TV)=D\times C^n$$ 
and 
$$\Phi _s((u|_{\partial D})^*TW)\subset S^1\times C^n$$
Moreover 
$$\Phi _0(u|_{\partial D})^*TW=S^1\times R^n$$
So, the homotopy induces 
a homotopy $\bar h$ in Lagrangian-Grassmanian
manifold. Note that $m(\bar h(0, \cdot ))=0$.
By the homotopy invariance of Maslov index,
we know that $m(u|_{\partial D})=0$.

\vskip 5pt 

   Consider the partial differential equation
\begin{eqnarray}
\bar \partial u+A(z)u=0  \ on \ D  \\
u(z)\in \Gamma (z) R^n\ for \ z\in \partial D \\
\Gamma (z)\in GL(2n,R)\cap Sp(2n)\\
m(\Gamma )=0 \ \ \ \ \ \ \ \ 
\end{eqnarray}

For $100<k<\infty $ consider the Banach space $\bar V_k $
consisting of all maps $u\in H^{k,2}(D, C^n)$ such 
that  $u(z)\in \Gamma (z)$ for almost all $z\in
\partial D$. Let $L_{k-1}$ the usual $L_{k-1}-$space $H_{k-1}(D,C^n)$ and 

$$L_{k-1}(S^1)=\{ u\in H^{k-1}(S^1)|u(z)\in \Gamma (z) R^n\ for \ z\in 
\partial D\}$$ 
We define an operator $P$:
$\bar V_{k}\rightarrow L_{k-1}\times L_{k-1}(S^1)$ by
\begin{equation}
P(u)=(\bar \partial u+Au,u|_{\partial D})
\end{equation}
where $D$ as in (2.1).
\begin{Proposition}
$\bar \partial : \bar V_p \rightarrow L_p$
is a real linear Fredholm operator of index n.
\end{Proposition}
Proof: see \cite{al,gro,wen}.

\section{Nonlinear Fredholm theory}

\subsection{Construction of Lagrangian submanifold}

Let $M$ be an open manifold and $(T^*M,p_idq_i)$ be 
the cotangent bundle of open manifold with the 
Liouville form $p_idq_i$. Since 
$M$ is open, there exists a function $g:M\to R$ without 
critical point. The translation by 
$tTdg$ along the fibre give a hamilton isotopy 
of $T^*M$:

\begin{equation}
h^T_t(q,p)=(q,p+tTdg(q))
\end{equation}

\begin{equation}
h^{T*}_t(p_idq_i)=p_idq_i+tTdg.
\end{equation}

\begin{Lemma}
For any given compact set $K\subset T^*M$, there exists 
$T=T_K$ such that $h^T_1(K)\cap K=\emptyset $.
\end{Lemma}
Proof. Similar to \cite{gro,ls}

\vskip 3pt 

Let $\Sigma \subset T^*M$ be a closed hypersurface, if there 
exists a vector field $V$ defined in the neighbourhood $U$ 
of $\Sigma $ transversal to $\Sigma $ 
such that $L_V\omega =\omega $, here 
$\omega =dp_i\wedge dq_i$ is a standard symplectic form  
on $T^*M$ induced by the Liouville form $p_idq_i$, we call 
$\Sigma $ the contact manifold of induced type in $T^*M$ 
with the induced contact form $\lambda =i_V\omega $. 

\vskip 3pt 

Let $(\Sigma ,\lambda )$ be a contact manifold of induced   
type or Weinstein's type in $T^*M$ with contact form 
$\lambda $ and $X$ its Reeb vector field, then 
$X$ integrates to a Reeb flow $\eta _s$ for $s\in R^1$.

By using the transversal vector field $V$, one 
can identify the neighbourhood $U$ of $\Sigma $ 
foliated by flow $f_t$ of $V$ and $\Sigma $, 
i.e., $U=\cup _tf_t(\Sigma )$ with the 
neighbourhood of $\{ 0\}\times \Sigma $ in the symplectization
$R\times \Sigma $ by the exact symplectic transformation(see\cite{ma2,vi1}).

Consider the form $d(e^a\lambda )$ 
at the point $(a,x)$ on the manifold 
$(R\times \Sigma )$, then one can check 
that $d(e^a\lambda )$ is a symplectic 
form on $R\times \Sigma $. Moreover 
One can check that 
\begin{eqnarray}
&&i_X(e^a\lambda )=e^a \\ 
&&i_X(d(e^a\lambda ))=-de^a 
\end{eqnarray}
So, the symplectization of Reeb vector field $X$ is the 
Hamilton vector field of $e^a$ with 
respect to the symplectic form $d(e^a\lambda )$. 
Therefore the Reeb flow lifts to the Hamilton flow 
$h_s$ on $R\times \Sigma $(see\cite{ag,eg}).

\vskip 3pt 

Let 
$$(V',\omega ')=(T^*M\times T^*M,d(p_i^1dq_i^1\ominus p_i^2dq_i^2))$$ 
be the anti-product of cotangent bundles and 
$${\cal L}=\{ (\sigma ,\sigma )|\sigma \in \Sigma \subset 
T^*M\}$$ 
be a closed 
isotropic submanifold 
contained in $
(\Sigma ',\lambda ')=(\Sigma \times \Sigma,\lambda \ominus \lambda )$, i.e.,  
there exists a smooth diagonal embedding $Q:{\cal L}\to \Sigma '$ such that 
$Q^*\lambda '|_{{\cal L}}=0$. 

Let
\begin{equation}
W'={\cal L}\times R, \ \ W'_s={\cal L}\times \{ s\} \label{eq:3.w}
\end{equation}
define 
\begin{eqnarray}
&&G':W'\to V'  \cr 
&&G'(w')=G'(l,s)=(\sigma ,\eta _s(\sigma )) \label{eq:3.ww}
\end{eqnarray}
we also denote $W'=G'(W')$

\begin{Lemma}
There does not exist any Reeb closed orbit in 
$(\Sigma ,\lambda )$ if and only if  
$G'(W'(s))\cap G'(W'(s'))$ is empty for $s\ne s'$.
\end{Lemma}
Proof. Obvious. 
\begin{Lemma}
If there does not exist any Reeb closed 
orbit for $X_\lambda $
in $(\Sigma ,\lambda )$ then 
there exists a smooth embedding 
$G':W'\to V'$ with $G'(l,s)=(\sigma ,\eta _s(\sigma ))$
such that
\begin{equation}
G'_K:{\cal L}\times (-K, K)\to V'  \label{eq:3.www}
\end{equation}
is a regular open Lagrangian embedding for any finite positive $K$.
\end{Lemma}
Proof. One first checks 
\begin{equation}
{G'}^*(e^a\lambda ')=\lambda _1-\eta (\cdot ,\cdot )^*\lambda _2
=\lambda _1-(\eta _s^*\lambda _2+i_{X}\lambda ds)=-ds
\end{equation}
Recall that $\Sigma $ is a contact manifold of induced type in $T^*M$, let 
$\lambda =i_Z(dp_i\wedge dq_i)$. Since 
$d\lambda =dp_i\wedge dq_i$, we know 
that 
$\theta =\lambda -p_idq_i$ is a close form  which determines  
a cohomology $[\theta ]\in H^1(\Sigma )$. 
Let 
$\theta _1=\lambda _1-p_i^1dq_i^1$
and $\theta _2=\lambda _2-p_i^2dq_i^2$. 
Since 
$[\theta _1]=[\theta _2]$, we have 
\begin{eqnarray}
(\lambda _1-p_i^1dq_i^1)-
(\lambda _2-p_i^2dq_i^2)=df
\end{eqnarray}
So,
\begin{eqnarray}
{G'}^*(p_i^1dq_i^1-p_i^2dq_i^2))=-ds-df.
\end{eqnarray}
This shows that $W'$ is an exact Lagrangian 
submanifold in $(T^*M\times T^*M,dp_i^1\wedge dq_i^1 \ominus 
dp_i^2\wedge dq_i^2)$.

Now we modify the above construction as follows\cite{mo}:  
\begin{eqnarray}
&&F_0':{\cal {L}}\times R\times R\to (R\times \Sigma )\times 
(R\times \Sigma )\cr 
&&F_0'(((0,\sigma ),(0,\sigma )),s,b)=((0,\sigma ),(b,\eta _s(\sigma )))
\end{eqnarray}
Now we embed a elliptic curve $E$ long along $s-axis$ and thin along $b-axis$ such that 
$E\subset [-s_1,s_2]\times [0,\varepsilon]$. We parametrize the $E$ by $t'$.

\begin{Lemma}
If there does not exist any closed Reeb orbit 
in $(\Sigma ,\lambda )$, 
then 
\begin{eqnarray}
&&F_0:{\cal {L}}\times S^1\to (R\times \Sigma )\times 
(R\times \Sigma )\cr 
&&F_0(((0,\sigma ),(0,\sigma )),t')=((0,\sigma ),(b(t'),\eta _{s(t')}(\sigma )))
\end{eqnarray}
is a compact Lagrangian submanifold. Moreover 
\begin{equation}
l(V',F_0({{\cal {L}}}\times S^1,d(e^a\lambda -e^b\lambda ))=area(E)
\end{equation}
\end{Lemma}
Proof. We check
that 
\begin{eqnarray}
{F_0}^*(e^a\lambda \ominus e^b\lambda )&=&-e^{b(t')}ds(t')
\end{eqnarray}
So, $F_0$ is a Lagrangian embedding.

If the circle $C$ homotopic to $C_1\subset {\cal {L}}\times s_0$ then  we compute
\begin{eqnarray}
\int _CF_0^*(e^b\lambda )=\int _{C_1}F_0^*(e^b\lambda )=0. 
\end{eqnarray}
since $\lambda |C_1=0$ due to $C_1\subset {\cal {L}}$ and 
$\cal L$ is Legendre submanifold. 

If the circle $C$ homotopic to $C_1\subset l_0\times S^1$ then  we compute
\begin{eqnarray}
\int _CF_0^*(e^b\lambda )=\int _{C_1}F_0^*(e^b\lambda )=n(area(E)). 
\end{eqnarray}
This proves the Lemma.

Now we modify the above construction as follows:  
\begin{eqnarray}
&&F':{\cal {L}}\times R\times R\to 
(0\times \Sigma )\times 
([0,\varepsilon ]\times \Sigma )\subset T^*M\times T^*M\cr 
&&F'(((0,\sigma ),(0,\sigma )),s,b)=((0,\sigma ),(b,\eta _s(\sigma )))
\end{eqnarray}
Now we embed a elliptic curve $E$ long along $s-axis$ and thin along $b-axis$ such that 
$E\subset [-s_1,s_2]\times [0,\varepsilon]$. We parametrize the $E$ by $t'$.

\begin{Lemma}
If there does not exist any closed Reeb orbit 
in $(\Sigma ,\lambda )$, 
then 
\begin{eqnarray}
&&F:{\cal {L}}\times S^1\to (R\times \Sigma )\times 
(R\times \Sigma )\cr 
&&F(((0,\sigma ),(0,\sigma )),t')=((0,\sigma ),(b(t'),\eta _{s(t')}(\sigma )))
\end{eqnarray}
is a compact Lagrangian submanifold. Moreover 
\begin{equation}
l(V',F({{\cal {L}}}\times S^1),d(p_i^1dq_i^1-p_i^2dq_i^2))=area(E)
\end{equation}
\end{Lemma}
Proof. We check
that 
\begin{eqnarray}
F^*(p_i^1dq_i^1-p_i^2dq_i^2)=
{F}^*(e^a\lambda \ominus e^b\lambda +df)=-e^{b(t')}ds(t')+df(l,t')
\end{eqnarray}
So, $F$ is a Lagrangian embedding.

If the circle $C$ homotopic to $C_1\subset {\cal {L}}\times s_0$ then  we compute
\begin{eqnarray}
\int _CF^*(p_i^1dq_i^1-p_i^2dq_i^2)=\int _{C_1}F^*(e^b\lambda +df)=0. 
\end{eqnarray}
since $\lambda |C_1=0$ due to $C_1\subset {\cal {L}}$ and 
$\cal L$ is ``Legendre'' submanifold. 

If the circle $C$ homotopic to $C_1\subset l_0\times S^1$ then  we compute
\begin{eqnarray}
\int _CF^*(p_i^1dq_i^1-p_i^2dq_i^2)=\int _{C_1}F^*(e^b\lambda +df)=n(area(E)). 
\end{eqnarray}
This proves the Lemma.

Now we construct an isotopy of Lagrangian embeddings as follows:  
\begin{eqnarray}
&&F':{\cal L}\times S^1\times [0,1]\to V'\cr 
&&F'(l,t',t)=(\sigma ,h^T_t(b(t'),\eta _{s(t')}(\sigma )))  \cr 
&&F'_t(l,t')=F'(l,t',t), \ l=(\sigma ,\sigma ).
\end{eqnarray}

\begin{Lemma}
If there does not exist any Reeb closed 
orbit for $X_\lambda $ 
in $(\Sigma ,\lambda )$ then $F'$ is an weakly exact isotopy of Lagrangian 
embeddings. Moreover 
for the choice of $T=T_\Sigma $ satisfying 
$\Sigma \cap h^T_1(\Sigma )=\emptyset $, then 
$F'_0({\cal L}\times S^1)\cap F'_1({\cal L}\times S^1)=\emptyset $.
\end{Lemma}
Proof. By Lemma3.1-3.5 and below.

\vskip 5pt

   Let 
$(V',\omega ')=(T^*M\times T^*M, dp_i^1\wedge dq_i^1
\ominus dp_i^2\wedge dq_i^2)$, 
$W'=F({\cal L}\times S^1)$, and 
$(V,\omega )=(V'\times C,\omega '\oplus \omega _0)$. 
As in \cite{gro}, we use symplectic figure eight trick invented by Gromov to 
construct a Lagrangian submanifold in $V$ through the 
Lagrange isotopy $F'$ in $V'$. 
Fix a positive $\delta <1$ and take a $C^{\infty }$-map $\rho :S^1\to 
[0,1]$, where the circle $S^1$ is parametrized by $\Theta \in [-1,1]$, 
such that the $\delta -$neighborhood $I_0$ of $0\in S^1$ goes to 
$0\in [0,1]$ and $\delta -$neighbourhood $I_1$ of $\pm 1\in S^1$ 
goes $1\in [0,1]$. Let 
\begin{eqnarray}
\tilde {l}&=&h^{T*}_{\rho }(p_i^1dq_i^1\ominus p_i^2dq_i^2)
=p_i^1dq_i^1-p_i^2dq_i^2-\rho (\Theta )Tdg\cr 
&=&-(e^{b(t')}ds(t')+d\beta )-\rho Tdg=(-e^{b(t')}ds(t')+d\beta +d\rho Tg)+Tgd\rho \cr 
&=&(-e^{b(t')}ds(t')+d\beta +d\rho Tg)-Tg\rho '(\Theta )d\Theta \cr
&=&(-e^{b(t')}ds(t')+d\beta +d\rho Tg)-\Phi d\Theta  
\end{eqnarray}
be the pull-back of the form 
$\tilde {l}'=(-e^{b(t')}ds(t')+d\beta +d\rho Tg)
-\psi (s,t)dt $ to $W'\times S^1$ under the map 
$(w',\Theta )\to (w',\rho (\Theta ))$ and 
assume without loss of generality $\Phi $
vanishes on $W'\times (I_0\cup I_1)$. Since 
$[\tilde {l}']|W'\times \{t\}=[(-e^{b(t')}ds(t')]$ is independent 
of $t$, so $F'$ is weakly exact. It is crucial here 
$|-\psi (s,t)|\leq M_0$  and $M_0$ is independent of $area(E)$.

  Next, consider a map $\alpha $ of the annulus $S^1\times [\Phi _-,\Phi _+]$ 
into $R^2$, where $\Phi _-$ and $\Phi _+$ are the lower and the upper 
bound of the function $\Phi $ correspondingly, such that 
   
   $(i)$ The pull-back under $\alpha $ of the form 
$dx\wedge dy$ on $R^2$ equals $-d\Phi \wedge d\Theta $. 
  
   $(ii)$ The map $\alpha $ is bijective on $I\times [\Phi _-,\Phi _+]$ 
where $I\subset S^1$ is some closed subset, 
such that $I\cup I_0\cup I_1=S^1$; furthermore, the origin 
$0\in R^2$ is a unique double point of the map $\alpha $ on 
$S^1\times 0$, that is 
$$0=\alpha (0,0)=\alpha (\pm 1,0),$$  
and 
$\alpha $ is injective on $S^1=S^1\times 0$ minus $\{ 0,\pm 1\}$. 

   $(iii)$ The curve $S^1_0=\alpha (S^1\times 0)\subset R^2$ ``bounds'' 
zero area in $R^2$, that is $\int _{S^1_0}xdy=0$, for the $1-$form 
$xdy$ on $R^2$. 
\begin{Proposition}
Let $V'$, $W'$ and $F'$ as above. Then there exists  
an exact Lagrangian embedding $F:W'\times S^1\to V'\times R^2$ 
given by $F(w',\Theta )=(F'(w',\rho (\Theta )),\alpha (\Theta ,\Phi ))$.
Denote $W=F(W'\times S^1)$. Then $W$ is contained in 
$T^*M\times T^*M\times B_R(0)$, here $4\pi R^2=8M_0$.
\end{Proposition}
Proof. Similar to \cite[2.3$B_3'$]{gro}.

\subsection{Formulation of Hilbert manifolds}

Let $(\Sigma ,\lambda )$ be a closed $(2n-1)-$ dimensional manifold
with a contact form $\lambda $ of induced type 
in $T^*M$, it is well-known that $T^*M$ is a Stein manifold, so 
it is exausted by a proper pluri-subharmonic function, in fact 
if $M$ is closed one can take $f={{1}\over {2}}|p|^2$, if $M$ 
is an open manifold one can take a proper Morse 
function $g$ to modify $f$, i.e., $f_1=f+\pi ^*g$. 
Since $\Sigma $ is compact and $W'=G'({\cal {L}}\times R)$ 
is contanied in $\Sigma $, by our construction 
we have $W$ is contained in a compact set $V_c$, 
$V_c\subset T^*M\times T^*M\times R^2$ for $M$ is an open manifold.  
If $M$ is closed and $\pi (\Sigma )=M$, we know 
that $W_K=F_K(W'_K\times 
S^1)$ is a bounded set in 
$T^*M\times T^*M\times R^2$. 

\vskip 3pt

We choose 
an almost complex structure $J_1$ on 
$T^*M$ tamed by $\omega _1 =dp_i\wedge dq_i$ and 
the metric $g_1=\omega _1(\cdot , J_1\cdot )$(see\cite{gro}).   
Let $(V',\omega ')=(T^*M\times T^*M,p_i^1dq_i^1\ominus p_i^2dp_i^2)$
By above discussion we know that  
$W'$ and  $\Sigma \times \Sigma $ 
contained in $\{ f_1\leq c\}\times \{ f_1\leq c\}$ 
for $c$ large 
enough, i.e., contained in a compact set $V'_c$ in $T^*M\times T^*M$. 
Then we expanding near $f_1^{-1}(c)$ to get 
a complete exact symplectic manifold with 
a complete Riemann metric with injective radius 
$r_0>0$(see\cite{ma2}).

\vskip 3pt 

In the following we denote by 
$(V,\omega )=(V'\times R^2,\omega '
\oplus dx\wedge dy))$ 
with the metric $g=g'\oplus g_0 $ induced by 
$\omega (\cdot ,J\cdot )$($J=J'\oplus i$ and 
$W\subset V$ a Lagrangian submanifold which was constructed in 
section 3.1.

   Let 
$${\cal D}^k(V,W,p)=\{ u \in H^k(D,V)|
u(x)\in W \ a.e \ for \ x\in \partial D \ and \ u(1)=p\}$$
for $k\geq 100$.
\begin{Lemma}
Let $W$ be a closed Lagrangian submanifold in 
$V$. Then, 
$${\cal D}^k(V,W,p)=\{ u \in H^k(D,V)|
u(x)\in W \ a.e \ for \ x\in \partial D \ and \ u(1)=p\}$$
is a pseudo-Hilbert manifold with the tangent bundle
\begin{equation}
T{\cal D}^k(V,W,p)=\bigcup _{u\in {\cal {D}}^k(V,W,p)}
\Lambda ^{k-1}(u^*TV,u|_{\partial D}^*TW,p)
\end{equation}
here 
$$\Lambda ^{k-1}(u^*TV,u|_{\partial D}^*TW,p)=$$ 
$$\{ H^{k-1}-sections \ of \ (u^*(TV),(u|_{\partial D})^*TL)\ 
which \ vanishes  \ at \ 1\} $$
\end{Lemma}
Proof: See \cite{al,kl}. 

\vskip 3pt

   Now we consider  a section
from ${\cal D}^k(V,W,p)$ to $T{\cal D}^k(V,W, p)$ follows as in 
\cite{al,gro}, i.e., 
let $\bar \partial :{\cal D}^k(V,W,p)\rightarrow T{\cal D}^k(V,W,p)$
be the Cauchy-Riemmann section 
\begin{equation}
\bar \partial u={{\partial u}\over {\partial s}}
+J{{\partial u}\over {\partial t}}  \label{eq:CR}
\end{equation}
for $u\in {\cal D}^k(V,W,p)$.

\begin{Theorem}
The Cauchy-Riemann section $\bar \partial $ defined in (\ref{eq:CR})
is a Fredholm section of Index zero.
\end{Theorem}
Proof. According to the definition of the Fredholm section, 
we need to prove that
$u\in {\cal D}^k(V,W,p)$, the linearization
$D\bar \partial (u)$ of $\bar \partial $ at $u$ is a linear Fredholm
operator.
Note that
\begin{equation}
D\bar \partial (u)=D{\bar \partial _{[u]}}
\end{equation}
where
\begin{equation}
(D\bar \partial _{[u]})v=\frac{\partial v}{\partial s}
+J\frac{\partial v}{\partial t}+A(u)v
\end{equation}
with 
$$v|_{\partial D}\in (u|_{\partial D})^*TW$$
here $A(u)$ is $2n\times 2n$
matrix induced by the torsion of
almost complex structure, see \cite{al,gro} for the computation.

   Observe that the linearization $D\bar \partial (u)$ of 
$\bar \partial $ at $u$ is equivalent to the following Lagrangian 
boundary value problem
\begin{eqnarray}
&&{{\partial v}\over {\partial s}}+J{{\partial v}\over {\partial t}}
+A(u)v=f, \ v\in \Lambda ^k(u^*TV)\cr 
&&v(t)\in T_{u(t)}W, \ \ t\in {\partial D}  \label{eq:Lin}
\end{eqnarray}
One 
can check that (\ref{eq:Lin}) 
defines a linear Fredholm operator. In fact, 
by proposition 2.2 and Lemma 2.1, since the operator $A(u)$ is a compact, 
we know that the operator $\bar \partial $ is a nonlinear Fredholm operator 
of the index zero.

\begin{Definition}
Let $X$ be a Banach manifold and $P:Y\to X$ the Banach 
vector bundle.
A Fredholm section $F:X\rightarrow Y$ is
proper if $F^{-1}(0)$ is a compact set and is called 
generic if $F$ intersects the zero section transversally, see \cite{al,gro}.
\end{Definition}
\begin{Definition}
$deg(F,y)=\sharp \{ F^{-1}(0)\} mod2$ is called the Fredholm
degree of a Fredholm section (see\cite{al,gro}).
\end{Definition}
\begin{Theorem}
Assum that the Fredholm section
$F=\bar \partial : {\cal D}^k(V,W,p)\rightarrow T^({\cal D}^k(V,W,p)$
constructed in (\ref{eq:CR}) is proper. Then,
$$deg(F,0)=1$$
\end{Theorem}
Proof: We assume that $u:D\mapsto V$ be a $J-$holomorphic disk
with boundary $u(\partial D)\subset W$ and 
by the assumption that $u$ is homotopic to the 
constant map $u_0(D)=p$. Since almost complex
structure ${J}$ tamed by  the symplectic form $\omega $,
by stokes formula, we conclude $u: D\rightarrow
V$ is a constant map. Because $u(1)=p$, We know that
$F^{-1}(0)={p}$.
Next we show that the linearizatioon $DF(p)$ of $F$ at $p$ is
an isomorphism from $T_p{\cal D}(V,W,p)$ to $E$.
This is equivalent to solve the equations
\begin{eqnarray}
{\frac {\partial v}{\partial s}}+J{\frac {\partial v}{\partial t}}
+Av=f\\
v|_{\partial D}\subset T_pW
\end{eqnarray}
here $J=J(p)=i$ and $A$ a constant matrix. By Lemma 2.1, we know that $DF(p)$ is an isomorphism.
Therefore $deg(F,0)=1$.

\section{Non-properness of a Fredholm section}

In this section we shall construct a non-proper Fredholm section 
$F_1:{\cal D}\rightarrow E$ by perturbing 
the Cauchy-Riemann section as in 
\cite{al,gro}.

\subsection{Anti-holomorphic section}

  Let $(V',\omega ')$ and $(V,\omega )=(V'\times C, 
\omega '\oplus \omega _0)$, and $W$ as in section3 
and $J=J'\oplus i$, $g=g'\oplus g_0$, 
$g_0$ the standard metric on $C$. 

   Now let $c\in C$ be a non-zero vector or nonzero 
constant vector field on $C$. We consider the 
equations
\begin{eqnarray}
v=(v',f):D\to V'\times C \nonumber \\
\bar \partial _{J'}v'=0,\bar \partial f=c\ \ \ \nonumber \\
v|_{\partial D}:\partial D\to W\ \ \  \label{eq:4.1}
\end{eqnarray}
here $v$ homotopic to constant map 
$\{ p\}$ relative to $W$. 
Note that $W\subset V\times B_R(0)$ for a positive 
number $R$ large enough. 
\begin{Lemma}
Let $v$ be the solutions of (\ref{eq:4.1}), then one has 
the following estimates
\begin{eqnarray}
E({v})=
\{ 
\int _D(g'({{\partial {v'}}\over {\partial x}},
{J'}{{\partial {v'}}\over {\partial x}})
+g'({{\partial {v'}}\over {\partial y}},
{J'}{{\partial {v'}}\over {\partial y}}) \nonumber \\
+g_0({{\partial {f}}\over {\partial x}},
{i}{{\partial {f}}\over {\partial x}})
+g_0({{\partial {f}}\over {\partial y}},
{i}{{\partial {f}}\over {\partial y}}))d\sigma \}
\leq 4\pi R^2. 
\end{eqnarray}
\end{Lemma}
Proof: Since $v(z)=(v'(z),f(z))$ satisfy (\ref{eq:4.1})
and $v(z)=(v'(z),f(z))\in V'\times C$ 
is homotopic to constant map $v_0:D\to \{ p\}\subset W$ 
in $(V,W)$, by the Stokes formula
\begin{equation}
\int _{D}v^*(\omega '\oplus \omega _0)=0
\end{equation}
Note that the metric $g$ is adapted to the symplectic form 
$\omega $ and $J$, i.e., 
\begin{equation}
g=\omega  (\cdot ,J\cdot )
\end{equation}
By the simple algebraic computation, we have 
\begin{equation}
\int _{D}{v}^*\omega  ={{1}\over {4}}
\int _{D^2}(|\partial v|^2 
-|\bar {\partial }v|^2)=0
\end{equation}
and 
\begin{equation}
|\nabla v|={{1}\over {2}}(
|\partial v|^2 +|\bar \partial v|^2 
\end{equation}
Then 
\begin{eqnarray}
E(v)&=&\int _{D} |\nabla v| \nonumber \\ 
      &=&\int _{D}\{ {{1}\over {2}}(
|\partial v|^2+|\bar \partial v|^2)\} d\sigma \nonumber \\ 
&=&\pi |c|_{g_0}^2
\end{eqnarray}
By the equations (\ref{eq:4.1}), 
one get 
\begin{equation}
\bar \partial f=c \ on \ D
\end{equation}
We have 
\begin{equation}
f(z)={{1}\over {2}}c\bar z+h(z)
\end{equation}
here $h(z)$ is a holomorphic function on $D$. Note that  
$f(z)$ is smooth up to the boundary $\partial D$, then, by 
Cauchy integral formula
\begin{eqnarray}
\int _{\partial D}f(z)dz&=&{{1}\over {2}}c\int _{\partial D}
\bar {z}dz+\int _{\partial D}h(z)dz \cr
&=&\pi ic
\end{eqnarray}
So, we have 
\begin{equation}
|c|={{1}\over {\pi}}|\int _{\partial D^2}f(z)dz|
\end{equation}
Therefore, 
\begin{eqnarray}
E(v)&\leq &\pi |c|^2
\leq {{1}\over {\pi }}|\int _{\partial D}f(z)dz|^2      \cr
&\leq &{{1}\over {\pi }}|\int _{\partial D}|f(z)||dz|^2   \cr
&\leq &4\pi |diam(pr_2(W))^2 \cr
&\leq &4\pi R^2.
\end{eqnarray}
This finishes the proof of Lemma.

\begin{Proposition}
For $|c|\geq 3R$, then the 
equations (\ref{eq:4.1})
has no solutions. 
\end{Proposition}
Proof. By (4.11), we have 
\begin{eqnarray} 
|c|&\leq &{{1}\over {\pi }}\int _{\partial D}|f(z)||dz|\cr
&\leq &{{1}\over {\pi }}\int _{\partial D}
diam(pr_2(W))||dz| \cr
&\leq &2R
\end{eqnarray}
It follows that $c=3R$ can not be obtained by 
any solutions.

\subsection{Modification of section $c$}

Note that the section $c$ is not a section of the 
Hilbert bundle in section 3 since $c$ is not 
tangent to the Lagrangian submanifold $W$, we must modify it as follows:

\vskip 3pt 

  Let $c$ as in section 4.1, we define 
\begin{eqnarray}
c_{\chi ,\delta }(z,v)=\left\{ \begin{array}{ll}
c \ \ \ &\mbox{if\  $|z|\leq 1-2\delta $,}\cr
0 \ \ \ &\mbox{otherwise}
\end{array}
\right. 
\end{eqnarray}
Then by using the cut off function $\varphi _h(z)$ and 
its convolution with section 
$c_{\chi ,\delta }$, we obtain a smooth section 
$c_\delta$ satisfying

\begin{eqnarray}
c_{\delta }(z,v)=\left\{ \begin{array}{ll}
c \ \ \ &\mbox{if\  $|z|\leq 1-3\delta $,}\cr
0 \ \ \ &\mbox{if\  $|z|\geq 1-\delta $.}
\end{array}
\right. 
\end{eqnarray}
for $h$ small enough, for the convolution theory see \cite{hor}.

   Now let $c\in C$ be a non-zero vector and 
$c_\delta $ the induced anti-holomorphic section. We consider the 
equations
\begin{eqnarray}
v=(v',f):D\to V'\times C \nonumber \\
\bar \partial _{J'}v'=0,\bar \partial f=c_\delta \ \ \ \nonumber \\
v|_{\partial D}:\partial D\to W\ \ \  \label{eq:4.16}
\end{eqnarray}
which is a slight modification of (\ref{eq:4.1})
Note that $W\subset V\times B_{R}(0)$. 
Then by repeating the same argument as section 4.1., we obtain 
\begin{Lemma}
Let $v$ be the solutions of (\ref{eq:4.16}) and $\delta $ 
small enough, then one has 
the following estimates
\begin{eqnarray}
E({v})\leq 4\pi R^2. 
\end{eqnarray}
\end{Lemma}
and

\begin{Proposition}
For $|c|\geq 3R$, then the 
equations (\ref{eq:4.16})
has no solutions. 
\end{Proposition}

\subsection{Modification of $J\oplus i$}

Let $(\Sigma ,\lambda )$ be a closed contact 
manifold
with a contact form $\lambda $ of 
induced type in $T^*M$. 
Let $J_M$ be an almost complex 
structure on $T^*M$ and 
$J_1=J_M\ominus J_M\oplus i$ the almost complex structure on 
$T^*M\times T^*M\times R^2$ tamed by 
$\omega '\oplus \omega _0$. 
Let $J_2$ be any almost complex structure on 
$T^*M\times T^*M\times R^2$.

\vskip 3pt 

Now we consider the almost conplex structure 
on the symplectic fibration $D\times V \to D$ which will 
be discussed in detail in section 5.1., see also 
\cite{gro}.

\begin{eqnarray}
J_{\chi ,\delta }(z,v)=\left\{ \begin{array}{ll}
i\oplus J _1\ \ \ &\mbox{if\  $|z|\leq 1-2\delta $,}\cr
i\oplus J_2 \ \ \ &\mbox{otherwise}
\end{array}
\right. 
\end{eqnarray}
Then by using the cut off function $\varphi _h(z)$ and 
its convolution with section 
$J_{\chi ,\delta }$, we obtain a smooth section 
$J_\delta$ satisfying

\begin{eqnarray}
J_{\delta }(z,v)=\left\{ \begin{array}{ll}
i\oplus J_1 \ \ \ &\mbox{if\  $|z|\leq 1-3\delta $,}\cr
i\oplus  {J_2} \ \ \ &\mbox{if\  $|z|\geq 1-\delta $.}
\end{array}
\right. 
\end{eqnarray}
as in section 
4.2.

Then as in section 4.2, one can also reformulation 
of the equations (\ref{eq:4.16}) and get similar 
estimates of Cauchy-Riemann equations, we leave it 
as exercises to reader.

\begin{Theorem}
The Fredholm sections $F_1=\bar \partial +c_\delta 
: {\cal  {D}}^k(V,W,p) \rightarrow T({\cal {D}}^k(V,W,p))$ is not proper 
for $|c|$ large enough.
\end{Theorem}
Proof. See \cite{al,gro}.

\section{$J-$holomorphic section}

Recall that $W\subset 
V=T^*M\times T^*M\times B_R(0)$ as in section 3. 
The Riemann metric $g$ on $V'\times R^{2}$ 
induces a metric $g|W$.

   Now let $c\in C$ be a non-zero vector and 
$c_\delta $ the induced anti-holomorphic section. We consider the 
nonlinear inhomogeneous equations (4.16) and 
transform it into $\bar J-$holomorphic map by 
considering its graph as in \cite{al,gro}.

Denote by $Y^{(1)}\to D\times V$ the bundle of homomorphisms $T_s(D)\to
T_v(V)$. If $D$ and $V$ are given the disk and the almost 
K\"ahler manifold, then
we distinguish the subbundle $X^{(1)}\subset Y^{(1)}$ which consists of
complex linear homomorphisms and we denote $\bar X^{(1)}\to D\times V$ the
quotient bundle $Y^{(1)}/X^{(1)}$. Now, we assign to each $C^1$-map $
v:D\to V$ the section $\bar \partial v$ of the bundle $\bar X^{(1)}$ over
the graph $\Gamma _v\subset D\times V$ by composing the differential of $v$
with the quotient homomorphism $Y^{(1)}\to \bar {X}^{(1)}$. If $c_\delta 
:D\times
V\to \bar X$ is a $H^k-$ section we write $\bar \partial v=c_\delta $ 
for the
equation $\bar \partial v=c_\delta |\Gamma _v$.

\begin{Lemma}
(Gromov\cite{gro})There exists a unique almost complex 
structure $J_g$ on $D\times V$(which 
also depends on the given structures in $D$ and in $V$), such that 
the (germs of) $J_\delta-$holomorphic sections $v:D\to D\times V$ are exactly and 
only the solutions 
of the equations $\bar \partial v=c_\delta $. Furthermore, the 
fibres $z\times V\subset D\times V$ are $J_\delta-$holomorphic(
i.e. the subbundles $T(z\times V)\subset T(D\times V)$ are $J_\delta-$complex) 
and the structure 
$J_\delta|z\times V$ equals the original structure on $V=z\times V$.
Moreover $J_\delta $ is tamed by $k\omega _0\oplus \omega $ for 
$k$ large enough which is independent of $\delta $.
\end{Lemma}

\section{Proof of Theorem 1.1}

\begin{Theorem}
There exists a non-constant $J-$holomorphic map $u: (D,\partial D)\to 
(V'\times C,W)$ with $E(u)\leq 4\pi R^2.$
\end{Theorem}
Proof.  
The results in 
section 4 shows the solutions of equations (4.16) must 
denegerate to a cusp curves, i.e., we obtain a Sacks-Uhlenbeck's 
bubble, i.e., $J-$holomorphic sphere or disk with boundary 
in $W$, the exactness of $T^*M\times T^*M\times R^2$ rules out 
the possibility of $J-$holomorphic sphere. So, we get a holomorphic disk. For the more detail, see the proof   
of Theorem 2.3.B in \cite{gro}.

\vskip 3pt 

{\bf Proof of Theorem 1.1}. By the assumption of Theorem 1.1, 
we know that the Lagrangian submanifold 
$W$ in $T^*M\times T^*M\times R^2$ is embedded. Moreover 
$l=l(T^*M\times T^*M\times R^2,W,\omega )=\inf 
\{\int _Df^*\omega >0|f:(D,\partial D)\to (T^*M\times T^*M\times R^2,W)\}=area(E)$. 
By Theorem6.1, $l\leq 4\pi R^2$. If $area(E)$ is large enough, this is a contradiction. 
This 
implies the assumption that ${\cal {L}}$ has no self-intersection 
point under Reeb flow does not hold.

\end{document}